\documentclass[11pt]{amsart}

\def\F{\mathbb{F}}
\def\Fp{\mathbb{F}_p}

\def\su{{\subseteq}}
\def\la{{\langle}}
\def\ra{{\rangle}}

\def\cFp{{\mathcal F}_p}

\def\w{{\omega}}
\def\gr{{\rm gr}}

\def\Proof{\noindent{\sl Proof.}\ }
\def\qed{{\hfill $\Box$ \medbreak}}
\usepackage{amssymb}

\newtheorem{defi}{Definition}[section]
\newtheorem{thm}[defi]{Theorem}
\newtheorem{lem}[defi]{Lemma}
\newtheorem{cor}[defi]{Corollary}

\newtheorem{prop}[defi]{Proposition}

\newtheorem*{thmm}{Theorem}

\begin{document}
\title{Isomorphism invariants of restricted enveloping algebras}

\author{\textsc{Hamid Usefi}}
\thanks{ The research is supported by an NSERC Postdoctoral Fellowship.}

\address{Department of Mathematics, University of British Columbia, 1984 Mathematics Road,
Vancouver, BC, Canada, V6T 1Z2}
\email{usefi@math.ubc.ca}

\begin{abstract}
Let $L$ and $H$ be finite-dimensional restricted Lie algebras over a perfect field $\F$
such that   $u(L)\cong u(H)$, where $u(L)$ is the  restricted  enveloping algebra
 of $L$. We prove that if $L$ is  $p$-nilpotent and  abelian, then $L\cong H$.
We deduce that  if $L$ is abelian and $\F$  is algebraically closed, then  $L\cong H$.
We use these results to prove the main result of this paper stating  that if $L$ is  $p$-nilpotent, then
 $L/L'^p+\gamma_3(L)\cong H/H'^p+\gamma_3(H)$.
\end{abstract}


  \maketitle

\section{Introduction}

Let $L$ be a restricted Lie algebra with the restricted enveloping algebra
$u(L)$. We shall say that a particular invariant of $L$ is
{\em determined} by $u(L)$, if every restricted Lie algebra $H$ also possesses this invariant whenever
$u(L)$ and $u(H)$ are isomorphic as associative algebras. In particular, the restricted isomorphism problem asks
whether the isomorphism type of $L$ is determined by $u(L)$. This problem is motivated by the classical isomorphism
problem for group rings: is every finite group $G$ determined by its integral group ring $\mathbb{Z} G$?
The survey article \cite{Sa} contains most of the development in this area.
In the late 1980's, Roggenkamp and  Scott \cite{RS87} and Weiss \cite{W} independently settled down the group ring problem for finite nilpotent groups.

There are close analogies between restricted Lie algebras and finite $p$-groups.
In particular, the restricted isomorphism problem is the Lie analogue of the modular
isomorphism problem that asks: given finite $p$-groups $G$ and $H$ with the property that $\Fp G\cong \Fp H$ can we deduce that $G\cong H$? Here, $\Fp$ denotes the field of $p$ elements.
There has been intensive investigation on the modular isomorphism problem, however
the main problem is rather far from being completely answered.
Unfortunately not every technique  from finite $p$-groups can be used for restricted Lie algebras. For example it is known that the class sums form a basis of the center of $\F G$. It then follows that the
center of $G$ is determined, see Theorem 6.6 in \cite{Seh}. Whether or not the center of $L$ is determined by $u(L)$ remains an interesting open question.

In analogy with finite $p$-groups we consider the class $\cFp$ of restricted Lie algebras that are finite-dimensional and $p$-nilpotent. Let $L\in \cFp$.  It follows from the Engel's Theorem that $L$ is nilpotent.
We shall examine the nilpotence class of $L$ in Corollary \ref{nilp-class}.
Note  that whether or not the nilpotence class of $G$ is determined by $\Fp G$ has been considered in the
recent years, however no major result is reported up-to-date, see \cite{BK}.

We start the investigation on the restricted isomorphism problem by first considering the abelian case. In
Proposition \ref{prop-perfect} we prove that if  $L\in \cFp$ is an abelian  restricted Lie algebra over a perfect
field  $\F$, then the isomorphism type of $L$ is determined by $u(L)$. Furthermore, if $\F$ is algebraically closed then every abelian restricted Lie algebra is determined by its enveloping algebra, see Corollary \ref{prop-alg-closed}.

It is not clear what is the next step beyond the abelian case in both the modular isomorphism
problem and the restricted isomorphism problem. Nevertheless, we have proved in \cite{U08} that if $L\in \cFp$ is a metacyclic restricted Lie algebra over a perfect field then the isomorphism type of  $L$ is determined by $u(L)$. The main result of this paper that will be proved in  Section \ref{quo-sec}, is another contribution in this direction; a similar result for finite $p$-groups was proved by Sandling \cite{S}. Let us recall that for a Lie subalgebra $I\su L$, we denote by  $I^p$ the restricted Lie subalgebra of $L$ generated by all $x^p$, $x\in I$. Also, $\gamma_i(L)$ denotes the $i$-th term of the lower central series of $L$. Our main result is as follows:

\begin{thmm} Let $L\in \cFp$ be a restricted Lie algebra over a perfect field. Then the restricted Lie algebra
$L/(L'^p+\gamma_3(L))$ is determined.
\end{thmm}

\section{Preliminaries}\label{pre-sec}
Let $L$ be a restricted Lie algebra with the restricted enveloping algebra $u(L)$ over a field $\F$.
By the Poincar\' e-Birkhoff-Witt (PBW) Theorem, see \cite{J}, we can view $L$ as a restricted Lie subalgebra of $u(L)$. Let $\w(L)$ denote the augmentation ideal of $u(L)$ which is the kernel of the augmentation map $\epsilon_{_L}: u(L)\to \F$  induced by $x\mapsto 0$, for every $x\in L$.

Let $H$ be another restricted Lie algebra such that
 $\varphi: u(L)\to u(H)$ is an algebra isomorphism.
We observe that the map $\eta: L\to u(H)$ defined by
$\eta=\varphi-\varepsilon _{_H}\varphi$ is a restricted Lie algebra homomorphism.
Hence,  $\eta$ extends to an algebra homomorphism
$\overline{\eta}: u(L)\to u(H)$. In fact,
$\overline{\eta}$ is an isomorphism that preserves the augmentation ideals, that is
$\overline{\eta}(\w(L))=\w(H)$, see \cite{RU} for the proof of similar fact for Lie algebras.
So, without loss of generality, we assume that $\varphi:u(L)\to u(H)$ is an
algebra isomorphism that preserves the augmentation ideals.

Recall that $L$ is said to be nilpotent if $\gamma_n(L)=0$ for
some $n$; the nilpotence class of $L$, denoted by $cl(L)$, is the minimal integer $c$
such that $\gamma_{c+1}(L)=0$. We denote by $L'_p$ the restricted subalgebra of $L$ generated by $L'=\gamma_2(L)$.
The $n$-th dimension subalgebra of $L$ is
$$
D_n(L)=L\cap \w^n(L)=\sum_{ip^j\geq n} \gamma_i(L)^{p^j},
$$
see \cite{RSh}.

Recall that $L$ is said to be in the class $\cFp$ if $L$ is finite-dimensional and $p$-nilpotent.
The \textit{exponent} of $x\in L$, denoted by $\exp(x)$, is the least integer $s$ such that
$x^{p^s}=0$. Whether or not $L\in \cFp$ is determined by the following lemma, see \cite{RSh}.
\begin{lem}\label{w(L)-nilpotent}
Let $L$ be a restricted Lie algebra. Then $L\in \cFp$ if and only if $\w(L)$ is nilpotent.
\end{lem}

 Now,  consider the graded restricted Lie algebra:
$$
\gr(L):=\bigoplus_{i\ge 1} D_i(L)/D_{i+1}(L),
$$
where the Lie bracket and the $p$-map are induced from $L$.
It is well-known that $u(\gr(L))\cong \gr(u(L))$ as algebras, see \cite{U08}.
So we may identify $\gr(L)$ as the graded restricted Lie subalgebra of $\gr(u(L))$
generated by $\w^1(L)/\w^2(L)$. Thus, $\gr(L)$ is determined. We can now deduce the following:
\begin{cor}\label{nilp-class}
 Let $L$ and $H$ be restricted Lie algebras such that $u(L)\cong u(H)$.
 If $L\in \cFp$ then $\mid cl(L)-cl(H)\mid \leq 1$.
 \end{cor}
 \Proof Let $c=cl(L)$.  We note that
$$
\gamma_n(\gr(L))=\bigoplus_{i\geq n}  \gamma_i(L){+}D_{i+1}(L)/D_{i+1}(L),
$$
for every $n\geq 1$. Since $\gr(L)$ is determined, it follows that $\gamma_{c+1}(\gr(H))=0$.
Hence, $\gamma_{c+1}(H)\su D_{c+2}(H)$. So, $\gamma_{c+2}(H)=\gamma_{c+3}(H)$.
Since $H$ is nilpotent, it follows that $\gamma_{c+2}(H)=0$.
\qed

Note that $
  D_n(\gr(L))=\bigoplus_{i\geq n}  D_i(L)/  D_{i+1}(L)$.
Thus, $  D_n(L)/ D_{n+1}(L)$ is determined, for every $n\geq 1$.  We remark that
methods of  \cite{RS} and \cite{RU}  can be adapted to prove that
 $D_n(L)/D_{2n+1}(L)$ and $D_{n}(L)/D_{n+2}(L)$ are also determined, for every $n\geq 1$.
In particular, $L/D_3(L)$ is determined.
We shall need the following  analogue of Lemma 5.1 in \cite{RU}.
\begin{lem}\label{Dn-wn}
If $\varphi:u(L)\to u(H)$ is an isomorphism then $\varphi(D_n(L)+\w^{n+1}(L))=D_n(H)+\w^{n+1}(H)$, for every
positive integer $n$.
\end{lem}

Now suppose that $L$ is an abelian restricted Lie algebra. Note that the conditions on  the $p$-map reduces to
$$
(x+y)^p=x^p+y^p, \quad (\alpha x)^p=\alpha^px^p,
$$
for every $x,y\in L$ and $\alpha\in \F$. Thus the $p$-map is a semi-linear transformation.
Let $\sigma$ be an automorphism of $\F$. Consider the skew polynomial ring $\F[t;\sigma]$ which consists of
polynomials $f(t)\in \F[t]$ with multiplication given by
$$
\alpha t^i \beta t^j= \alpha\beta^{\sigma^{-i}} t^{i+j}.
$$
It is well-known that $\F[t;\sigma]$ is a PID.
Now suppose that  $\F$ is perfect and let
$\sigma$ be the automorphism given by
 $\sigma(\alpha)=\alpha^p$.
Note that  $\F[t; \sigma]$ acts on $L$ by $ x\cdot t= x^p$. So, by the theory of finitely generated modules over a PID,
 $L$ decomposes as a direct sum of cyclic $\F[t; \sigma]$-modules. In particular, the number of these summands is unique.
 We summarize this in the following, see also \cite{J} or Section 4.3 in \cite{BMPZ}. We denote by $\langle x\rangle_{p}$
 the subalgebra generated by $x$.

\begin{thm}\label{abelian-rest-perfect}
Let $L$ be a finitely generated abelian restricted Lie algebra over a perfect field $\F$. Then there exist
a unique integer $n$ and generators $x_1,\ldots,x_n\in L$ such that
$$
L=\langle x_1\rangle_{p}\oplus\cdots \oplus \langle x_n\rangle_{p}.
$$
\end{thm}

\begin{prop}\label{prop-perfect}
Let $L\in \cFp$ be an abelian  restricted Lie algebra over a perfect field  $\F$. If $H$ is  a restricted Lie
algebra such
that $u(L)\cong u(H)$, then $L\cong H$.
\end{prop}
\Proof
We argue by induction on $\text{dim}_{\F} L$. Let $A$ be the
subalgebra of $\w(L)$ generated by all $u^p$, where $u\in \w(L)$.
We observe that $A\cong \w(L^p)$,  as algebras. Thus there is an induced isomorphism:
$$
\w(L^p)\cong \w(H^p).
$$
Since $L\in \cFp$, it follows that  $\text{dim}_{\F} L^p<\text{dim}_{\F} L$.
Thus, by the induction hypothesis, there exists a restricted Lie algebra isomorphism $\phi: L^p\cong H^p$. We
now lift $\phi$ to an isomorphism of $L$ and $H$.
By Theorem \ref{abelian-rest-perfect}, there exist  generators $x_1,\ldots, x_n\in L$ such that
$L=\langle x_1\rangle_{p}\oplus\cdots \oplus \langle x_n\rangle_{p}.$ Without loss of generality we assume
$$
L^p=\langle x_1^p\rangle_{p}\oplus\cdots \oplus \langle x_m^p\rangle_{p},
$$
for some $m\leq n$. Thus, $x_i^p=0$, for every $i$ in the range $m<i\leq n$.
Note that $\dim L=n+\dim L^p$. So, as it is mentioned in Theorem \ref{abelian-rest-perfect}, $n$ is determined.
Let $y_1,\ldots, y_n\in H$ such that
$H=\langle y_1\rangle_{p}\oplus\cdots \oplus \langle y_n\rangle_{p}.$ Then
$$
H^p=\langle y_1^p\rangle_{p}\oplus\cdots \oplus \langle y_m^p\rangle_{p}.
$$
So, we can assume that  $\phi(x_i^p)=y_i^p$, for every $1\leq i\leq m$. We can verify that the map induced by $x_i\mapsto y_i$, for every $1\leq i\leq n$, is a restricted Lie algebra isomorphism between $L$ and $H$.
\qed

\begin{cor}\label{L/L'_p}
Let  $L\in \cFp$ be a restricted Lie algebra over a perfect field. Then $L/L'_p$ is  determined.
\end{cor}
\Proof Note that $[u(L),u(L)]u(L)=L'_pu(L)$. Also, we have
$u(L/L'_p)\cong u(L)/L'_pu(L)$. Hence, $u(L/L'_p)$ is determined.
Since $L/L'_p\in \cFp$, it follows from Proposition \ref{prop-perfect} that $L/L'_p$ is determined.
\qed

It turns out that over an  algebraically closed field stronger results hold. Before we state the next result we need to recall a well-known theorem, see \cite{J} or Section 4.3 in \cite{BMPZ}. Let  $T_{_L}= \langle x\in L \mid x^p=x\rangle_{\F}$ and denote  by $ \mathrm{Rad}(L)$  the subalgebra of $L$ spanned by all $p$-nilpotent elements.

\begin{thm}\label{abelian-rest-algclosed}
Let $L$ be a finite-dimensional abelian restricted Lie algebra over an algebraically closed  field $\F$. Then
$L=T_{_L}\oplus \mathrm{Rad}(L).$
\end{thm}

\begin{cor}\label{prop-alg-closed}
Let $L$ be a finite-dimensional abelian  restricted Lie algebra over an algebraically closed field $\F$. Let $H$
be a restricted Lie algebra such that $u(L)\cong u(H)$. Then $L\cong H$.
\end{cor}
 \Proof
Note that for every $k\geq 1$,
$$
\dim_{\F} L/D_{p^k}(L)=\dim_{\F} L/D_{p}(L)+\cdots+\dim_{\F} D_{p^{k-1}}(L)/D_{p^k}(L),
$$
is determined. So
$\dim_{\F}D_{p^k}(L)$ is determined, for every $k\geq 1$. Let $t$ be the least integer such that
$\mathrm{Rad}(L)^{p^t}=0$.
It follows that $D_{p^t}(L)=T_{_L}$ . Hence,
$\dim_{\F}\mathrm{Rad}(L)=\dim_{\F}\mathrm{Rad}(H)$, by Theorem \ref{abelian-rest-algclosed}.
  Note that $L/T_{_L}\cong \mathrm{Rad}(L)$, as restricted Lie algebras.
 We claim that $\varphi(u(T_{_L}))=u(T_{_H})$. Suppose that the claim holds. Then $\varphi(T_{_L}u(L))=T_{_H}u(H)$. So,
\begin{align*}
 u(L/T_{_L})\cong u(L)/T_{_L}u(L)\cong u(H)/T_{_H}u(H)\cong u(H/T_{_H}).
\end{align*}
Thus, $u(\mathrm{Rad}(L))\cong u(\mathrm{Rad}(H))$.
Since $\mathrm{Rad}(L),\mathrm{Rad}(H)\in \cFp$, Proposition \ref{prop-perfect} then implies that there exists an
  isomorphism $\phi: \mathrm{Rad}(L)\to \mathrm{Rad}(H)$. Clearly, $\phi$ can be extended to an
isomorphism of $L$ and $H$ by sending $x_i$ to $y_i$.

Now, we prove the claim. Let $z_1,\ldots,z_n$ be a basis of $\mathrm{Rad}(H)$ and $y_1,\ldots,y_s$ be a basis of $T_{_H}$
 and assume that every $y_i$ is less than every $z_j$. Let $x\in
T_{_L}$ and express $\varphi(x)$ in terms of PBW monomials in the $y_i$ and $z_j$.
So we have,
$$
\varphi(x)=u+\sum \alpha y_1^{a_1}\ldots y_s^{a_s}z_1^{b_1}\ldots z_n^{b_n},
$$
 where $u$ is a linear combination of  PBW monomials in the $y_i$ only and each term in the   sum has the
property that $b_1+\cdots+b_n\neq 0$. Note that for a large $k$ we have
$\varphi(x)^{p^k}=u^{p^k}\in u(T_{_H}).$
But $\varphi(x)=\varphi(x)^{p^k}$. So, $\varphi(x)\in u(T_{_H})$.
 Since $u(T_{_L})$ is generated by $L$ and $\varphi$ is an algebra homomorphism, it follows that
$\varphi(u(T_{_L}))\su u(T_{_H})$. But $u(T_{_L})$ and $u(T_{_H})$ are finite-dimensional. So we get
$\varphi(u(T_{_L}))=
u(T_{_H})$.  This proves the claim and so the proof is complete.\qed

\section{The quotient $L/L'^p+\gamma_3(L)$}\label{quo-sec}

We first record a couple of easy statements.

\begin{lem}
Let $N$ be a restricted subalgebra of $L$. We have,
$$
\w(L)N+N\w(L)=[N,L]+N\w(L)
$$
\end{lem}

\begin{lem}\label{Nu(L)}
For every restricted subalgebra $N$ of $L$ we have,
\begin{enumerate}
\item $L\cap ([N,L]+N\w(L))=[N,L]+N^p$.
\item $Nu(L)/\w(L)N+N\w(L)\cong N/([N,L]+N^p)$.
\end{enumerate}
\end{lem}
Now write $J_L=\w(L)L'+L'\w(L)=\w(L)L'_p+L'_p\w(L)$. Since both $\w(L)L'$ and $L'\w(L)$ are determined, it
follows that $J_L$ is determined.

\begin{cor}\label{dim-L-mod}
If $L\in \cFp$ then $\dim_{\F} (L/L'^p{+}\gamma_3(L))$ is determined.
\end{cor}
\Proof Since $L'_pu(L)$ and $J_L$ are determined, it follows from  Lemma \ref{Nu(L)} that
 $\dim_{\F} (L'_p/L'^p{+}\gamma_3(L))$ is determined. The result then follows, since  $L/L'_p$ is determined,
by Corollary \ref{L/L'_p}. \qed

From now on we assume that $L\in\cFp$ and $\F$ is perfect.
By Theorem \ref{abelian-rest-perfect}, there exists $e_1, \ldots, e_n\in L$ such that
$$
L/L'_p=\langle  e_1+L'_p\rangle_{p}\oplus\cdots \oplus \langle e_n+L'_p\rangle_{p}.
$$
 Let $\bar X$ be a basis of $L/L'_p$ consisting of $\bar e_i^{p^j}$, where $\bar e_i=e_i+L'_p$ and
$1\leq i\leq n$. We fix a set $X$ of representatives of $\bar X$. So the elements of $X$ are
linearly independent modulo $L'_p$.

We define the \textit{height} of an element $x\in L$, denoted by $\nu(x)$,
to be  the largest integer $n$ such that $x\in D_n(L)$, if $n$ exists and infinity otherwise.
The \textit{weight} of a PBW monomial $x_{1}^{a_1}\ldots x_{t}^{a_t}$ is defined to be $\sum_{i=1}^t a_i \nu
(x_i)$.  We observe that  $\nu(e_i^{p^j})=p^{j}$, for every $1\leq i\leq n$ and every $1\leq j<\exp (\bar e_i)$.
Indeed, if $e_i^{p^j}\in D_{m}(L)$, for some $m> p^{j}$,
then $e_i^{p^j}=\sum_{k>j} \alpha_k e_i^{p^k}$ modulo $L'_p$.
It follows then that  $e_i^{p^{\exp(\bar e_i)-1}}\in L'_p$, which is a contradiction.
Let $Y$ be a  linearly independent subset of $L'_p$  such that
$Z=X\cup Y$ is a basis of $L$ and the  set
 $\{z+D_{\nu(z)+1}\mid z\in Z\}$ is a basis of $\gr (L)$.
 One way to construct such a subset $Y$ is to take coset representatives of a basis for
  $$
 \oplus_{i\geq 1} D_i(L)\cap(L'_p+\la X\ra_{\F})/D_{i+1}(L).
 $$
We need  the following variant of Theorem 2.1 in \cite{RSh}.

\begin{lem}\label{lem-RSh}
Let $L\in \cFp$. Let $\bar Z$ be a homogeneous basis of $\gr(L)$ with a fixed set of representatives $Z$.
Then the set of all PBW monomials in $Z$ of weight at least $k$  forms a basis for $\w^k(L)$, for every $k\geq
1$.
\end{lem}

Note that $J_L$ is linearly independent with the set of all PPW monomials in $X$.
Let $E$ denote the vector space spanned by $J_L$ and all PBW monomials
in $X$ of degree at least two.
The following lemma is easy to see and so we omit the proof.
\begin{lem}\label{w(L)-decom}
The following statements hold.
\begin{enumerate}
\item $\w(L)=L+E$.
\item $(L+J_L)\cap E=J_L=E\cap L'_p u(L)$.
\item $\w(L)/J_L= L+J_L/J_L\oplus E/J_L$.
\end{enumerate}
\end{lem}

\begin{lem}\label{E/J_L}
If $L\in \cFp$ then  $E/J_L$ is a central restricted Lie ideal of $\w(L)/J_L$.
\end{lem}
\Proof The fact that $E/J_L$ is a central  Lie ideal of $\w(L)/J_L$ easily follows from the identity
$[ab,c]=a[b,c]+[a,c]b$ which holds in any associative algebra.
So we have to prove that $E/J_L$ is closed under the $p$-map. Since $J_L$ is an associative ideal of
$\w(L)$,  it is enough to prove that $u^p\in E$, for every PBW monomial $u$ in $E$.
Let  $u=e_1^{a_1}\cdots e_n^{a_n}$, where each $a_i$ is in the range $0\leq a_i < p^{\exp{\bar e_i}}$.
 It is not hard to see that
$u^p=e_1^{pa_1}\cdots e_n^{pa_n}$ modulo $J_L$.
Since $L\in \cFp$, each $\bar e_i$ is $p$-nilpotent.
If $pa_i< p^{\exp(\bar e_i)}$, for every $1\leq i\leq n$, then $u^p$
is a PBW monomial of degree at least two.
Now suppose that $pa_i\geq p^{\exp(\bar e_i)}$, for some $i$.
If $pa_i=p^{\exp(\bar e_i)}$ then $a_i$ is a power of $p$.
Since $u$  has degree at least two, there exists $j\neq i$ such that $a_j\neq 0$.
It now follows that $u^p\in J_L$.
 If $pa_i>p^{\exp(\bar e_i)}$ then $e_i^{pa_i}\in J_L$ and so  $u^p\in E$.
\qed

\begin{lem}\label{H-phi(E)}
We have $H\cap \varphi(E)\su J_H$.
\end{lem}
\Proof We suppose  $J_H$=0 and prove that $H\cap \varphi(E)=0$. Let $v\in H\cap \varphi(E)\su \w^2(H)$.
Let $u\in E$ such that $\varphi(u)=v$. So, $u\in \w^2(L)$. We prove by induction that $u\in \w^{n}(L)$, for every $n$.
But $\w(L)$ is nilpotent, by Lemma \ref{w(L)-nilpotent}, and so $u=0$.
Suppose now, by induction, that $u\in w^{n}(L)$
and we prove that $u\in w^{n+1}(L)$.
So, $v\in H\cap w^{n}(H)=D_n(H)$. Thus, by Lemma \ref{Dn-wn}, $u\in (D_{n}(L)+\w^{n+1}(L))\cap E$.
But
$$
(D_n(L)+\w^{n+1}(L))\cap E\su \w^{n+1}(L).
$$
Indeed, let $u=\sum \alpha _i z_i+w$, where each $z_i\in Z$ has height $n$ and $w\in \w^{n+1}(L)$.
By Lemma \ref{lem-RSh}, $w$ is a linear combination of PBW  monomials  in $Z$ of weight at least $n+1$.
Since $u\in E$ it follows by the PBW Theorem that $\alpha_i=0$, for every  $i$. So, $u=w\in \w^{n+1}(L)$, as required.
\qed

\begin{lem}
We have,  $\w(H)/J_H= H{+}J_H/J_H\oplus \varphi(E)/J_H$.
\end{lem}
\Proof
By Lemma \ref{H-phi(E)}, it is enough to prove that
 $$
 \w(H)/J_H\su H+J_H/J_H\oplus \varphi(E)/J_H.
 $$
 Note that both $\w(H)/J_H$ and $\varphi(E)/J_H$ are determined. Since $\dim_{\F} (H+J_H/J_H)=\dim_{\F}
(H/(H')^p+\gamma_3(H))$ is determined by Corollary \ref{dim-L-mod}, the
 result follows from Lemma \ref{w(L)-decom}.
 \qed

We can now finish the proof of our main result. Note that $L+J_L/J_L\cong L/L'^p+\gamma_3(L)$,
by Lemma \ref{Nu(L)}.
\begin{lem}
The restriction of the natural isomorphism $\w(L)/J_L\to \w(H)/J_H$ to $L+J_L/J_L$ induces an
isomorphism of $L+J_L/J_L$ and $H+J_H/J_H$.
\end{lem}
\Proof
We denote by $\varphi$ the induced isomorphism  $\w(L)/J_L \to \w(H)/J_H$.
Let $\varphi_{\mid_{L+J_L/J_L}}=\varphi_1+\varphi_2$ denote the restriction of $\varphi$ to $L+J_L/J_L$,
where
$\varphi_1: L+J_L/J_L  \to H+J_H/J_H$.
It is  enough to
 prove that $\varphi_1$  is a restricted Lie algebra isomorphism.
 Since $E/L$ is a central  Lie ideal of $\w(L)/J_L$, by lemma \ref{E/J_L},
$\varphi(E)/J_H$ is a central  Lie ideal of $\w(H)/J_H$. So, for every  $x,z\in L$, we have
$$
\varphi([x,z]+J_L)=[\varphi(x)+J_H,\varphi(z)+J_H]=[\varphi_1(x),\varphi_1(z)]+J_H.
$$
So, $\varphi_1$ preserves the Lie brackets. Also,
\begin{align*}
\varphi(x^p+J_L)&=\varphi(x)^p+J_H
=(\varphi_1(x))^p+(\varphi_2(x))^p+J_H
\end{align*}
Since $(\varphi_2(x))^p+J_H\in\varphi(E)/J_H$, it follows that $\varphi_1$ preserves the $p$-powers.
Furthermore, $\varphi_1$ is injective, by Lemma \ref{w(L)-decom}. Since $L+J_L/J_L$ and
$H+J_H/J_H$ have the same dimension, by Corollary \ref{dim-L-mod}, it follows that $\varphi_1$ is an isomorphism,
as required.
\qed

\section*{Acknowledgments}
I am grateful to the referee for careful reading of the paper and Luzius Grunenfelder for useful discussions.

\end{document}